\documentclass[a4paper, 11pt]{amsart}
\usepackage{amsfonts, amsthm, amssymb, amsmath, stmaryrd}
\usepackage{mathrsfs,array}
\usepackage{eucal,color,times,enumerate,accents}
\usepackage{tikz-cd}
\usepackage{url}
\usepackage{bbm}
\usepackage[utf8]{inputenc}
\usetikzlibrary{decorations.text}
\usepackage{pstricks}
\usetikzlibrary{decorations.pathmorphing, patterns,shapes}
\usetikzlibrary{arrows,shapes,positioning}
\usetikzlibrary{decorations.markings}
\tikzstyle arrowstyle=[scale=1]
\tikzstyle directed=[postaction={decorate,decoration={markings,
    mark=at position .65 with {\arrow[arrowstyle]{stealth}}}}]
\tikzstyle reverse directed=[postaction={decorate,decoration={markings,
    mark=at position .65 with {\arrowreversed[arrowstyle]{stealth};}}}]
\usepackage{pgfplots}
\usepackage{lipsum}
\usepackage{graphicx}
\usepackage{rotating}
\usepackage{enumitem}
\usepackage{enumerate}
\usepackage{textcomp}
\usepackage[pdftex,bookmarksnumbered,bookmarksopen]{hyperref}
\hypersetup{colorlinks, linkcolor=blue}
\usepackage[colorinlistoftodos]{todonotes}
\usepackage[left=4cm,right=4cm,top=3cm,bottom=3cm]{geometry}
\usepackage{mathtools}
\input xy
\xyoption{all}
\usepackage{yfonts}
\usepackage{hyphenat}
\newcommand{\F}{\mathbb{F}}

\newcommand{\Q}{\mathbb{Q}}

\newcommand{\Ql}{\mathbb{Q}_{\ell}}

\newcommand{\Image}{\mathrm{Im}}
\newcommand{\Coker}{\mathrm{Coker}}

\newcommand{\etale}{\'etale }

\begin{document}
\bibliographystyle{alpha}
\newtheorem{theorem}[equation]{Theorem}
\newtheorem{hope}[equation]{Hope}
\newtheorem{proposition}[equation]{Proposition}
\newtheorem{lemma}[equation]{Lemma}
\newtheorem{claim}[equation]{Claim}
\newtheorem{corollary}[equation]{Corollary}

\newtheorem{fact}[equation]{Fact}

\theoremstyle{definition}
\newtheorem{definition}[equation]{Definition}
\newtheorem{question}[equation]{Question}
\newtheorem{conjecture}[equation]{Conjecture}
\newtheorem{answer}[equation]{Answer}
\newtheorem{remark}[equation]{Remark}
\newtheorem{example}[equation]{Example}
\newtheorem{warning}[equation]{Warning}
\newtheorem{notation}[equation]{Notation}
\newtheorem{construction}[equation]{Construction}

\title{A note on the behaviour of the Tate conjecture under finitely generated field extensions}
\author{Emiliano Ambrosi}

\begin{abstract}
We show that the $\ell$-adic Tate conjecture for divisors on smooth proper varieties over finitely generated fields of positive characteristic follows from the $\ell$-adic Tate conjecture for divisors on smooth projective surfaces over finite fields. Similar results for cycles of higher codimension are given.
\end{abstract}
\maketitle
\numberwithin{equation}{subsection} 

\section{Introduction}
Let $k$ be a field of characteristic $p\geq 0$ with algebraic closure $\overline k$ and write $\pi_1(k)$ for the absolute Galois group of $k$. A $k$-variety is a reduced scheme, separated and of finite type over $k$. For a $k$-variety $Z$ write $Z_{\overline k}:=Z\times_k\overline k$ and $CH^i(Z_{\overline k})$ for the group of algebraic cycles of codimension $i$ modulo rational equivalence. Let $\ell\neq p$ be a prime. 
\subsection{Conjectures}
Recall the following versions of the Grothendieck-Serre-Tate conjectures (\cite{tateconj}, \cite[Section 7.3]{andremotivi}):
\begin{conjecture}\label{aminconjtate}
If $k$ is finitely generated and $Z$ is a smooth proper $k$-variety, then:
\begin{itemize}
\item $T(Z,i,\ell):$ The $\ell$-adic cycle class map $$c_{Z_{\overline k}}: CH^i({Z_{\overline k}})\otimes \Ql\rightarrow \bigcup_{[k':k]<+\infty}H^{2i}(Z_{\overline k},\Ql(i))^{\pi_1(k')}$$ is surjective;
\item $S(Z,i,\ell):$ The action of $\pi_1(k)$ on $H^{2i}(Z_{\overline k},\Ql(i))$ is semisimple;
\item $WS(Z,i,\ell):$ The inclusion $H^{2i}(Z_{\overline k},\Ql(i))^{\pi_1(k)}\subseteq H^{2i}(Z_{\overline k},\Ql(i))$ admits a $\pi_1(k)$-equivariant splitting.
\end{itemize}
\end{conjecture}
For a field $K$, one says that $T(K,i,\ell)$ holds if for every finite field extension $K\subseteq L$ and every smooth proper $L$-variety $Z$, $T(Z,i,\ell)$ holds. One defines similarly the conditions $S(K,i,\ell)$ and $WS(K,i,\ell)$. 

Conjecture \ref{aminconjtate} is widely open in general. By the works of many people, $T(Z,1,\ell)$ is known when $Z$ is an abelian variety (\cite{tateabelianfinite}, \cite{zarhinabelian1}, \cite{zarhinabelian2}, \cite{faltingsabelian}), a K3 surface (\cite{Nygaardk3}, \cite{tank3}, \cite{andk3}, \cite{charlietatek3}, \cite{peratatek3}, \cite{peratatek32}) or when $Z$ lies in some other special class of $k$-varieties; see for example \cite[Section 5.13]{peratatek3} and \cite{Moonentate}. 
\subsection{Behaviour under finitely generated field extensions}
\subsubsection{Main result}
For abelian varieties and K3 surfaces, Conjecture \ref{aminconjtate} is closely related to the finiteness of rational points on their moduli spaces; see \cite[Proposition 2]{tateabelianfinite} and \cite{finitenessk3}. This may suggest that Conjecture \ref{aminconjtate} could be easier to prove when $k$ is a finite field. The main result of this note is that, to prove Conjecture \ref{aminconjtate} for varieties over finitely generated fields of positive characteristic, it is actually enough to prove it for varieties over finite fields.
\begin{theorem}\label{higher}
If $p>0$, then $T(\F_p,i,\ell)$ and $WS(\F_p,i,\ell)$ imply $T(k,i,\ell)$ for every finitely generated field $k$ of characteristic $p$.
\end{theorem}

\subsubsection{The case of divisors}
By \cite[Proposition 2.6.]{Tate}, if algebraic and numerical equivalences on algebraic cycles coincide rationally in degree $i$, then $T(Z,i,\ell)$ implies $WS(Z,i,\ell)$. Since this holds for divisors, from Theorem \ref{higher} we deduce the following.
\begin{corollary}\label{reductionTate}
If $p>0$, then $T(\F_p,1,\ell)$ implies $T(k,1,\ell)$  for every finitely generated field $k$ of characteristic $p$.
\end{corollary}
By an unpublished result (\cite{djtate}) of De Jong (whose proof has been simplified in \cite[Theorem 4.3]{Varmor}), over finite fields the $\ell$-adic Tate conjecture for divisors for smooth projective varieties follows from the $\ell$-adic Tate conjecture for divisors for smooth projective surfaces. Hence Corollary \ref{reductionTate} implies the following:
\begin{corollary}\label{reductionTate2}
If $p>0$, then $T(Z,1,\ell)$ for every finite field $\F_q$ of characteristic $p$ and every smooth projective $\F_q$-surface $Z$ implies $T(k,1,\ell)$ for every finitely generated field $k$ of characteristic $p$.
\end{corollary}
\subsection{Previously known results}
We quickly review previously known results on the behaviour of Conjecture \ref{aminconjtate} under finitely generated field extension. Let $k\subseteq K$ be a regular finitely generated field extension and let $Z$ be a smooth proper $K$-variety. 
Choose a geometrically connected, smooth $k$-variety $\mathcal K$ with generic point $\eta:K\rightarrow \mathcal K$ and a smooth proper morphism $\mathfrak f:\mathcal Z\rightarrow \mathcal K$ fitting into a cartesian diagram:
\begin{center}
\begin{tikzcd}
Z\arrow{d}\arrow{r}{i_{\eta}}\arrow[phantom]{rd}{\Box}&  \mathcal Z\arrow{d}{\mathfrak f} \\
k\arrow{r}{\eta} & \mathcal K.
\end{tikzcd}
\end{center}
For every $\mathfrak t\in \mathcal K$, choose a geometric point $\overline{\mathfrak t}$ over $\mathfrak t$ and write $\mathcal Z_{\mathfrak t}$ and $\mathcal Z_{\overline{\mathfrak t}}$ for the fibre of $\mathfrak f$ at $\mathfrak t$ and $\overline{\mathfrak t}$ respectively.
\subsubsection{Characteristic zero}\label{andreiiii}
Assume that $k$ has characteristic zero. Then the arguments in \cite[Section 5.1]{AndreIHES}, show that $S(k,i,\ell)$ and $T(k,i,\ell)$ imply $T(K,i,\ell)$. The idea is use resolution of singularities to embed $\mathcal Z$ as a dense open subset into a smooth proper $k$-variety $\mathcal Z^{cmp}$. Then, by smooth proper base change the action of $\pi_1(K)$ on $H^{2i}(Z_{\overline K},\Ql(i))$ factors through the surjection $\pi_1(K)\twoheadrightarrow \pi_1(\mathcal K)$. By the global invariant cycles theorem (\cite{Weil2}; see \cite[Theoreme 1.1.1]{Andre}), the natural map $H^{2i}(\mathcal Z_{\overline k}^{cmp},\Ql(i))\rightarrow H^{2i}(Z_{\overline K},\Ql(i))^{\pi_1(\mathcal K_{\overline k})}$ is surjective. Since $S(\mathcal Z_{\overline k}^{cmp},i,\ell)$ holds, the map $$H^{2i}(\mathcal Z_{\overline k}^{cmp},\Ql(i))^{\pi_1(k)}\rightarrow H^{2i}(Z_{\overline K},\Ql(i))^{\pi_1(K)}$$ is still surjective. Hence, by $T(\mathcal Z^{cmp},i,\ell)$, every class in $H^{2i}(Z_{\overline K},\Ql(i))^{\pi_1(K)}$ arises from a cycle in $H^{2i}(\mathcal Z_{\overline k}^{cmp},\Ql(i))$.

Let us point out that, even assuming $S(\F_p,i,\ell)$ instead of $WS(\F_p,i,\ell)$ in Theorem \ref{higher}, the arguments in \cite[Section 5.1]{AndreIHES} do not work in positive characteristic, since resolution of singularities is not known. On the other hand our arguments for Theorem \ref{higher} do not work in characteristic zero, since they use in an essential way the procylicity of $\pi_1(\F_q)$.
\subsubsection{Semisimplicity}
As observed in \cite{leifusiccomeimmobile}, Deligne's geometric semisimplicity theorem \cite[Theoreme 3.4.1]{Weil2} can be used to show that $S(k,i,\ell)$ implies $S(K,i,\ell)$. More precisely, by Deligne's geometric semisimplicity theorem (\cite[Theoreme 3.4.1]{Weil2}), the restriction of the action of $\pi_1(\mathcal K)$ on $H^{2i}(Z_{\overline K},\Ql(i))$ to its normal subgroup $\pi_1(\mathcal K_{\overline k})\subseteq \pi_1(\mathcal K)$ is semisimple. For any closed point $\mathfrak t\in \mathcal K$, the subgroup of $\pi_1(\mathcal K)$ generated by $\pi_1(\mathcal K_{\overline k})$ and by the image of $\pi_1(\mathfrak t)\rightarrow \pi_1(\mathcal K)$ is open in $\pi_1(\mathcal K)$. So, since the action of $\pi_1(\mathcal K_{\overline k})$ is semisimple, the action of $\pi_1(\mathcal K)$ on $H^{2i}(Z_{\overline K},\Ql(i))$ is semisimple if the action of $\pi_1(\mathfrak t)$ on $H^{2i}(Z_{\overline K},\Ql(i))$ induced via restriction through the morphism $\pi_1(\mathfrak t)\rightarrow \pi_1(\mathcal K)$ is semisimple. But this action identifies, modulo the isomorphism $H^{2i}(Z_{\overline K},\Ql(i))\simeq H^{2i}(\mathcal Z_{\overline{\mathfrak t}},\Ql(i))$ given by the choice of an \etale path between $\overline \eta:\overline K\rightarrow \mathcal K$ and $\overline{\mathfrak t}$, with the natural action of $\pi_1(\mathfrak t)$ on $H^{2i}(\mathcal Z_{\overline{\mathfrak t}},\Ql(i))$. Hence $S(\mathcal Z_{\mathfrak t},i,\ell)$ implies $S(Z,i,\ell)$.
\subsubsection{Infinite finitely generated fields}
Assume now that $k$ is \textbf{infinite} and finitely generated. Then the results in \cite{AndreIHES} (see \cite[Corollary 5.4]{MC}) if $p=0$ or \cite[Theorem 1.3.3]{mioneron} if $p>0$, show that $T(k,1,\ell)$ implies $T(K,1,\ell)$. Indeed, they show that there exists always a closed fibre $\mathcal Z_{\mathfrak t}$ such that the N\'eron-Severi group $NS(\mathcal Z_{\mathfrak t})\otimes \Q$ of $\mathcal Z_{\mathfrak t}$ identifies (rationally) with the N\'eron-Severi group $NS(Z)\otimes \Q$ of $Z$. Since the choice of an \etale path between $\overline{\mathfrak t}$ and $\overline \eta$ induces a commutative diagram of injective maps
\begin{center}
\begin{tikzcd}
NS(Z)\otimes \Ql\arrow[hook]{r}{\simeq}\arrow[hook]{d}&NS(\mathcal Z_{\mathfrak t})\otimes \Ql\arrow[hook]{d}\\
H^2(Z_{\overline K},\Ql)^{\pi_1(K)}\arrow[hook]{r} &H^2(Z_{\overline{\mathfrak t}},\Ql(1))^{\pi_1(\mathfrak t)},
\end{tikzcd}
\end{center}
this shows that $T(\mathcal Z_{\mathfrak t},1,\ell)$ implies $T(Z,1,\ell)$. 
\subsection{Acknowledgements}
The author is very grateful to a referee for pointing out gap in a previous version of this note and for his accurate comments, which helped clarify the exposition of the paper. This paper is part of the author’s Ph.D. thesis under the supervision of Anna Cadoret. He is grateful for her careful re-readings of this paper and her constructive suggestions. He thanks Marco D'Addezio for helpful discussions around the use of \cite[Proposition 2.6.]{Tate}.
\section{Proof of Theorem \texorpdfstring{\ref{higher}}-}\label{moving}
Let $k$ be an infinite finitely generated field $k$ of characteristic $p>0$ inside a fixed algebraic closure $\overline k$. Let $\F_q$ (resp. $\F$) the algebraic closure of $\F_p$ in $k$ (resp. $\overline k$). For every smooth proper $k$-variety $Z$, write $CH^{i}_{\ell}(Z_{\overline k})$ for the image of $CH^{i}(Z_{\overline k})\otimes \Ql\rightarrow H^{2i}(Z_{\overline k},\Ql(i))$. Let $Z$ be a smooth and proper $k$-variety.
\subsection{Strategy}\label{sectionstrategyTate}
Compared to the arguments Section in \ref{andreiiii}, the extra difficulties one has to deal with in the proof of Theorem \ref{higher} come from the fact that resolution of singularities is not known in positive characteristic and that we don't assume the semisimplicity of the Galois action in $\ell$-adic cohomology.

To overcome the use of resolution of singularities, we use De Jong's alterations theorem to construct a generically \etale alteration $\widetilde{\mathcal Z}\rightarrow \mathcal Z$ such that $\widetilde{\mathcal Z}$ embeds as a dense open subset into a smooth proper $\F_q$-variety. As pointed out by a referee, the resulting morphism $\widetilde{\mathcal Z}\rightarrow \mathcal Z\rightarrow \mathcal K$ is not, in general, generically smooth, so that we cannot apply directly the global invariant cycles theorem. To solve this issue, we use the main ingredients of its proof: the Hard Lefschetz theorem \cite[Theorem 4.1.1]{Weil2} and the theory of weights for $\F_q$-varieties \cite[Theorem 1]{Weil2}.

To overcome the lack of the semisimplicity assumption we combine the procyclicty of $\pi_1(\F_q)$ with the condition $WS(Z,i,\ell)$, to study the fixed points of the action of $\pi_1(\F_q)$ (Section \ref{usingtate}) via the generalized eigenspace of generalized eigenvalue $1$ of a topological generator.
\subsection{Preliminary reductions}\label{sprel}
To prove $T(Z,i,\ell)$, one may freely replace $k$ with a finite field extension. In particular we may assume that all the connected components of $Z_{\overline k}$ are defined over $k$ and so, working with each component separately, that $Z$ is geometrically connected over $k$. The following well known lemma, a slight variant of \cite[Theorem 5.2]{Tate}, will be used twice.
\begin{lemma}\label{alteration}
Let $W$ be a smooth proper $k$-variety and $g:W\rightarrow Z$ a generically finite dominant morphism. Then the following hold:
\begin{itemize}
\item The map $g^*:H^{2i}(Z_{\overline k},\Ql(i))\rightarrow H^{2i}(W_{\overline k},\Ql(i))$ is injective.
\item For any $z\in H^{2i}(Z_{\overline k},\Ql(i))$, if $g^*(z)\in CH^{i}_{\ell}(W_{\overline k})$ then $z\in CH^{i}_{\ell}(Z_{\overline k})$.
\end{itemize}  
In particular $T(W,i,\ell)$ implies $T(Z,i,\ell)$.
\end{lemma}
\proof
Assume first that $W$ is geometrically connected. Then, by Poincar\'e duality, there is a morphism $g_*:H^{2i}(W_{\overline k},\Ql(i))\rightarrow H^{2i}(Z_{\overline k},\Ql(i))$ which is compatible with the push forward of cycles $g_*:CH^i(W_{\overline k})\otimes \Ql\rightarrow CH^i(Z_{\overline k})\otimes \Ql$ and such that $g_*g^*$ is equal to the multiplication by the generic degree of $g:W\rightarrow Z$. All the assertions then follow from the commutative diagram:
\begin{center}
\begin{tikzcd}
CH^i(Z_{\overline k})\otimes \Ql\arrow{r}{g^*}\arrow{d}{c_{Z_{\overline k}}}& CH^i(W_{\overline k})\otimes \Ql\arrow{r}{g_*}\arrow{d}{c_{W_{\overline k}}}& CH^i(Z_{\overline k})\otimes \Ql\arrow{d}{c_{Z_{\overline k}}}\\
H^{2i}(Z_{\overline k},\Ql(i))\arrow{r}{g^*}& H^{2i}(W_{\overline k},\Ql(i))\arrow{r}{g_*}& H^{2i}(Z_{\overline k},\Ql(i)).
\end{tikzcd}
\end{center}
In general, we reduce to the situation where $W$ is geometrically connected. To prove Lemma \ref{alteration}, we can freely replace $k$ with a finite field extension and hence assume that all the connected components  $W_{i,\overline k}$ of $W_{\overline k}$ are defined over $k$. Since $g:W\rightarrow Z$ is dominant and generically finite and $Z$ is connected, there is at least one connected component (say $W_1$) mapping surjectively onto $Z$. Since $Z$ and $W_1$ are smooth proper $k$-varieties of the same dimension, the morphism $g_1:W_1\rightarrow W\rightarrow Z$ is still dominant and generically finite. The general case follows then from the geometrically connected case and the diagram:
\begin{center}
\begin{tikzcd}
CH^i(Z_{\overline k})\otimes \Ql\arrow{r}\arrow{d}&CH^i(W_{\overline k})\otimes \Ql\arrow{r}\arrow{d}&CH^i(W_{1,\overline k})\otimes \Ql\arrow{d}\\
H^{2i}(Z_{\overline k},\Ql(i))\arrow{r}&H^{2i}(W_{\overline k},\Ql(i))\arrow{r}&H^{2i}(W_{1,\overline k},\Ql(i)).\qedhere
\end{tikzcd}
\end{center}
By De Jong's alteration theorem (\cite{alteration}) applied to $Z_{\overline k}$, there exists a smooth projective $\overline k$-variety $W'$ and a dominant generically finite morphism $g':W'\rightarrow Z_{\overline k}$. By descent and replacing $k$ with a finite field extension, there exist a smooth projective $k$-variety $W$ and a dominant generically finite morphism $g:W\rightarrow Z$ which, after base change along $Spec(\overline k)\rightarrow Spec(k)$, identifies with $g':W'\rightarrow Z_{\overline k}$. By Lemma \ref{alteration}, we may replace $Z$ with $W$ and hence we may assume that $Z$ is a smooth projective $k$-variety.
Since the action of $\pi_1(k)$ on $CH^{i}_{\ell}(Z_{\overline k})$ factors through a finite quotient, replacing $k$ with a finite field extension, we may and do assume that $CH^{i}_{\ell}(Z_{\overline k})\subseteq H^{2i}(Z_{\overline k},\Ql(i))^{\pi_1(k)}$. The core of the proof is the following proposition.
\begin{proposition}\label{propositionkey}
Let $Z$ be a geometrically connected smooth projective $k$-variety such that $CH^{i}_{\ell}(Z_{\overline k})\subseteq H^{2i}(Z_{\overline k},\Ql(i))^{\pi_1(k)}$. Assume that $T(\F_p,i,\ell)$ holds. Replacing $k$ with a finite field extension, there exist a projective $k$-scheme $\widetilde Z$ and a dominant generically finite morphism $h:\widetilde Z\rightarrow Z$, such that for every $z\in H^{2i}(Z_{\overline k},\Ql(i))^{\pi_1(k)}$ one has $h^*(z)\in CH^{i}_{\ell}(\widetilde Z_{\overline k})$.
\end{proposition}
Before proving Proposition \ref{propositionkey}, let us show that it implies Theorem \ref{higher}. Replacing $k$ with a finite field extension we can take $h:\widetilde Z\rightarrow Z$ as in the statement of Proposition \ref{propositionkey}. Write $\widetilde Z_{\overline k, red}$ for the reduced closed subscheme of $\widetilde Z_{\overline k}$. Then $h_{red}:\widetilde Z_{\overline k,red}\rightarrow \widetilde Z_{\overline k} \rightarrow Z_{\overline k}$ is still dominant and generically finite and for every $z\in H^{2i}(\widetilde Z_{\overline k},\Ql(i))^{\pi_1(k)}$ one has $h_{red}^*(z)\in CH^{i}_{\ell}(\widetilde Z_{\overline k,red})$. So, by descent and replacing $k$ with a finite field extension, we can assume that $\widetilde Z$ is geometrically reduced and that all the irreducible components of $\widetilde Z_{\overline k}$ are defined over $k$. Then, by De Jong alteration's theorem applied to $\widetilde Z_{\overline k}$ and descent, replacing $k$ with a finite field extension, there exists a generically finite dominant morphism $W\rightarrow \widetilde Z$ with $W$ a smooth projective $k$-variety. The morphism $g:W\rightarrow \widetilde Z\rightarrow Z$ is still generically finite and dominant and for every $z\in H^{2i}(Z_{\overline k},\Ql(i))^{\pi_1(k)}$ one has $g^*(z)\in CH^{i}_{\ell}(W_{\overline k})$. Then Theorem \ref{higher} follows from Lemma \ref{alteration}. \endproof
The next subsection is devoted to the proof of Proposition \ref{propositionkey}.
\subsection{Proof of Proposition \ref{propositionkey}}
\numberwithin{equation}{subsubsection} 
Let $Z$ be a geometrically connected smooth projective $k$-variety such that $CH^{i}_{\ell}(Z_{\overline k})\subseteq H^{2i}(Z_{\overline k},\Ql(i))^{\pi_1(k)}$.
\subsubsection{Spreading out and alterations}\label{sspread}
Spreading out to $\F_q$, there exist a geometrically connected, smooth $\F_q$-variety $\mathcal K$ with generic point $\eta:k\rightarrow \mathcal K$ and a smooth projective morphism $\mathfrak f:\mathcal Z\rightarrow \mathcal K$ with geometrically connected fibres fitting into a cartesian diagram:
\begin{center}
\begin{tikzcd}
Z\arrow{d}\arrow{r}{i_{\eta}}\arrow[phantom]{rd}{\Box}&  \mathcal Z\arrow{d}{\mathfrak f} \\
k\arrow{r}{\eta} & \mathcal K.
\end{tikzcd}
\end{center}
By De Jong alteration's theorem, there exist an integral smooth $\F_q$-variety $\widetilde{\mathcal Z}$, an open embedding $\widetilde{\mathfrak i}:\widetilde{\mathcal Z}\rightarrow \widetilde{\mathcal Z}^{cmp}$ with dense image into a smooth projective $\F_q$-variety $\widetilde{\mathcal Z}^{cmp}$ and a generically \'etale, proper, dominant morphism $\mathfrak h:\widetilde{\mathcal Z}\rightarrow \mathcal Z$. Then $\widetilde{\mathcal Z}^{cmp}$ is geometrically connected over some finite field extension $\F_q\subseteq \F_{q'}$. Replacing $\F_q$ with $\F_{q'}$ amounts to replacing $k$ with the finite field extension $k':=k\F_{q'}$, so we can assume that $\widetilde{\mathcal Z}$ and $\widetilde{\mathcal Z}^{cmp}$ are geometrically connected over $\F_q$.

Since $\widetilde{\mathcal Z}\rightarrow \mathcal Z\rightarrow \F_q$ is quasi-projective, the morphism $\mathfrak h:\widetilde{\mathcal Z}\rightarrow \mathcal Z$ is quasi-projective as well (\cite[Tag 0C4N]{stacks-project}). Since $\mathfrak f:\mathcal Z\rightarrow \mathcal K$ is projective, this implies that $\widetilde{\mathcal Z}\rightarrow \mathcal K$ is quasi-projective. Since $\mathfrak h:\widetilde{\mathcal Z}\rightarrow \mathcal Z$ and $\mathfrak f:\mathcal Z\rightarrow \mathcal K$ are proper, the morphism $\widetilde{\mathcal Z}\rightarrow \mathcal K$ is proper as well. So $\widetilde{\mathcal Z}\rightarrow \mathcal K$ is proper and quasi-projective hence projective. 
The generic fibre $\widetilde Z\rightarrow k$ of $\widetilde{\mathcal Z}\rightarrow \mathcal K$ is then a projective $k$-scheme endowed with a generically finite dominant morphism $h:\widetilde Z\rightarrow Z$. The situation is summarized in the following diagram of $\F_q$-schemes:
\begin{center}
\begin{tikzcd}
\widetilde{Z}\arrow{d}{h}\arrow{r}{\widetilde{i}_{\eta}}\arrow[phantom]{rd}{\Box} &  \widetilde{\mathcal Z}\arrow[hook]{r}{\widetilde{\mathfrak i}}\arrow{d}{\mathfrak h} & \widetilde{\mathcal Z}^{cmp}\\
Z\arrow{d}\arrow{r}{i_{\eta}}\arrow[phantom]{rd}{\Box}&  \mathcal Z\arrow{d}{\mathfrak f} \\
k\arrow{r}{\eta} & \mathcal K.
\end{tikzcd}
\end{center}
The edge map 
$$\mathrm{Ler}:H^{2i}(\mathcal Z_{\F},\Ql(i))\rightarrow H^0(\mathcal K_{\F},R^{2i}\mathfrak f_{\F*}\Ql(i))$$
in the Leray spectral sequence for $\mathfrak f_{\F}:\mathcal Z_{\F}\rightarrow \mathcal K_{\F}$ fits then into a commutative diagram:
\begin{center}
\begin{tikzcd}[column sep=small]
 & & &CH^i(\widetilde Z_{\overline k})\otimes \Ql\arrow[swap, near start]{d}{c_{\widetilde Z_{\overline k}}} &CH^i(Z_{\overline k})\otimes \Ql\arrow[swap]{l}{h^*}\arrow{d}{c_{Z_{\overline k}}} \\

CH^i(\widetilde{\mathcal Z}_{\F}^{cmp})\otimes \Ql\arrow{r}{\widetilde{\mathfrak i}^*}\arrow{d}{c_{\widetilde{\mathcal Z}_{\F}^{cmp}}}&CH^i(\widetilde{\mathcal Z}_{\F})\otimes \Ql\arrow{d}{c_{\widetilde{\mathcal Z}_{\F}}}\arrow{rru}{\widetilde i^*_{\eta}}&CH^i(\mathcal Z_{\F})\otimes \Ql \arrow{l}{\mathfrak h^*}\arrow[swap, near start]{d}{c_{\mathcal Z_{\F}}}\arrow[near end, swap]{urr}{i^*_{\eta}}&  H^{2i}(\widetilde Z_{\overline k},\Ql(i)) & H^{2i}(Z_{\overline k},\Ql(i))\arrow[swap]{l}{h^*}\\

H^{2i}(\widetilde{\mathcal Z}^{cmp}_{\F},\Ql(i))\arrow{r}{\widetilde{\mathfrak i}^*}&H^{2i}(\widetilde{\mathcal Z}_{\F},\Ql(i))\arrow[near end, swap]{urr}{\widetilde i^*_{\eta}}&H^{2i}(\mathcal Z_{\F},\Ql(i))\arrow[swap]{urr}{i^*_{\eta}}\arrow{l}{\mathfrak h^*}\arrow{rr}{\mathrm{Ler}} && H^0(\mathcal K_{\F},R^{2i}\mathfrak f_{\F*}\Ql(i))\arrow[hook]{u}.
\end{tikzcd}
\end{center}
\subsubsection{Hard Lefschetz Theorem}\label{sspread2}
Write $\varphi\in \pi_1(\F_q)$ for the arithmetic Frobenius of $\F_q$ and, for every $\pi_1(\F_q)$-representation $V$, write $V^{\varphi}_{gen}$ for the generalized eigenspace on which $\varphi$ acts with generalized eigenvalue $1$. Since $\pi_1(\F_q)$ is procyclic the $\pi_1(\F_q)$-equivariant inclusion $V^{\varphi}_{gen}\subseteq V$ as a $\pi_1(\F_q)$-equivariant splitting. Hence if $r:V\rightarrow W$ is $\pi_1(\F_q)$-equivariant morphism of continuous $\pi_1(\F_q)$-representations, one has 
\begin{equation}\label{usoprociclico}
\Image(r:V^{\varphi}_{gen}\rightarrow W)=\Image(r:V\rightarrow W)\cap W^{\varphi}_{gen}.
\end{equation}

Let $z$ be in $H^{2i}(Z_{\overline k},\Ql(i))^{\pi_1(k)}$. In this section we lift $h^*(z)\in H^{2i}(\widetilde{Z}_{\F},\Ql(i))$ to  $H^{2i}(\widetilde{\mathcal Z}_{\F},\Ql(i))^{\varphi}_{gen}$.
By smooth proper base change, the action of $\pi_1(k)$ on $H^{2i}(Z_{\overline k},\Ql(i))$ factors through the canonical surjection $\pi_1(\mathcal K)\rightarrow \pi_1(k)$, hence $H^{2i}(Z_{\overline k},\Ql(i))^{\pi_1(k)}\simeq H^{2i}(Z_{\overline k},\Ql(i))^{\pi_1(\mathcal K)}$. 
Consider the diagram:
\begin{center}
\begin{tikzcd}
H^{2i}(\widetilde{\mathcal Z}^{cmp}_{\F},\Ql(i))\arrow{r}{\widetilde{\mathfrak i}^*}&H^{2i}(\widetilde{\mathcal Z}_{\F},\Ql(i))&H^{2i}(\mathcal Z_{\F},\Ql(i))\arrow{r}{\mathrm{Ler}}\arrow[swap]{l}{\mathfrak h^*} & H^0(\mathcal K_{\F},R^{2i}\mathfrak f_{\F*}\Ql(i)).
\end{tikzcd}
\end{center}
Since $$z\in H^{2i}(Z_{\overline k},\Ql(i))^{\pi_1(k)}\simeq H^0(\mathcal K_{\F},R^{2i}\mathfrak f_{\F*}\Ql(i))^{\pi_1(\F_q)}\subseteq H^0(\mathcal K_{\F},R^{2i}\mathfrak f_{\F*}\Ql(i)),$$
the element $z$ is in $H^0(\mathcal K_{\F},R^{2i}\mathfrak f_{\F*}\Ql(i))^{\varphi}$.
Recall the following consequence of the Hard Lefschetz Theorem:
\begin{fact}\label{hardlefschetz}
The map 
$$\mathrm{Ler}: H^{2i}(\mathcal Z_{\F},\Ql(i))\rightarrow H^0(\mathcal K_{\F},R^{2i}\mathfrak f_{\F*}\Ql(i))$$
is surjective
\end{fact}
\proof
Let $\mathcal F\in Pic(\mathcal Z_{\F})$ be a relative ample line bundle for $\mathfrak f_{\F}:\mathcal Z_{\F}\rightarrow \mathcal K_{\F}$, write $m$ for the relative dimension of $\mathfrak f_{\F}:\mathcal Z_{\F}\rightarrow \mathcal K_{\F}$ and, to simplify the notation, choose an identification $\Ql(i)\simeq \Ql$. 
For every integer $r\geq 1$, $c_{\mathcal Z_{\F}}(\mathcal F)\in H^2(\mathcal Z_{\F},\Ql)$ induces, by cup product, a morphism
$$u^r:R^{m-r}\mathfrak f_{\F*}\Ql\rightarrow R^{m+r}\mathfrak f_{\F*}\Ql.$$
By proper base change and the Hard Lefschetz Theorem \cite[Theorem 4.1.1]{Weil2}, $u^r:R^{m-r}\mathfrak f_{\F*}\Ql\rightarrow R^{m+r}\mathfrak f_{\F*}\Ql$  is an isomorphism. In the terminology of \cite[(2.2)]{Spectraldeligne}, this shows that $\Ql$ satisfies the Lefschetz condition with respect to $c_{\mathcal Z_{\F}}(\mathcal F)$. By \cite[Proposition 2.4]{Spectraldeligne}, this implies that the Leray spectral sequence for $\mathfrak f_{\F}:\mathcal Z_{\F}\rightarrow \mathcal K_{\F}$
$$E_2^{a,b}:=H^{a}(\mathcal Z_{\F},R^{b}\mathfrak f_{\F*}\Ql)\Rightarrow H^{a+b}(\mathcal K_{\F},\Ql)$$
degenerates at the second page. Hence the edge map $H^{2i}(\mathcal Z_{\F},\Ql)\rightarrow H^0(\mathcal K_{\F},R^{2i}\mathfrak f_{\F*}\Ql)$ is surjective.\endproof
By Fact \ref{hardlefschetz}, the map $\mathrm{Ler}:H^{2i}(\mathcal Z_{\F},\Ql(i))\rightarrow H^0(\mathcal K_{\F},R^{2i}\mathfrak f_{\F*}\Ql(i))$ is surjective, hence, by  (\ref{usoprociclico}), $z\in H^0(\mathcal K,R^{2i}\mathfrak f_*\Ql(i))\subseteq H^0(\mathcal K_{\F},R^{2i}\mathfrak f_{\F*}\Ql(i))^{\varphi}_{gen}$ is the image of some $z'\in H^{2i}(\mathcal Z_{\F},\Ql(i))^{\varphi}_{gen}$ and then $\mathfrak h^*(z')\in H^{2i}(\widetilde{\mathcal Z}_{\F},\Ql(i))^{\varphi}_{gen}$. 
\subsubsection{Theory of weights}
We now prove that $\mathfrak h^*(z')$ is the image of some $\widetilde z\in H^{2i}(\widetilde{\mathcal Z}^{cmp}_{\F},\Ql(i))^{\varphi}_{gen}$ via $\widetilde{\mathfrak i}^*:H^{2i}(\widetilde{\mathcal Z}^{cmp}_{\F},\Ql(i))\rightarrow H^{2i}(\widetilde{\mathcal Z}_{\F},\Ql(i))$. Write $d$ for the common dimension of $\widetilde{\mathcal Z}$ and $\widetilde{\mathcal Z}^{cmp}$. The localization exact sequence for the dense open immersion $\widetilde{\mathcal Z}\rightarrow \widetilde{\mathcal Z}^{cmp}$ with complement $\mathcal D:=\widetilde{\mathcal Z}^{cmp}-\widetilde{\mathcal Z}$, gives an exact sequence
$$H^{2d-2i-1}_c(\mathcal D_{\F},\Ql(-i))(d)\rightarrow H^{2d-2i}_c(\widetilde{\mathcal Z}_{\F},\Ql(-i))(d) \rightarrow H^{2d-2i}_c(\widetilde{\mathcal Z}_{\F}^{cmp},\Ql(-i))(d).$$
Combining this sequence with Poincar\'e duality for the smooth varieties $\widetilde{\mathcal Z}$ and $\widetilde{\mathcal Z}^{cmp}$, one sees that
\begin{equation}\label{inclusionequation}
\Coker\Big(\text{ }\widetilde{\mathfrak i}^*:H^{2i}(\widetilde{\mathcal Z}^{cmp}_{\F},\Ql(i))\rightarrow H^{2i}(\mathcal Z_{\F},\Ql(i))\text{ }\Big)\subseteq (H^{2d-2i-1}_c(\mathcal D_{\F},\Ql(-i))(d))^{\vee}.
\end{equation}
We combine (\ref{inclusionequation}) with the theory of weights.
\begin{fact}[{\cite[Th\'eor\`eme 3.3.1 and Corollaire 3.3.9]{Weil2}}]\label{weights}
Let $\mathcal X$ be a separated scheme of finite type over $\F_q$. Then, for every integers $m\geq 0$ and $n$, $H_c^m(\mathcal X_{\F},\Ql(n))$ is mixed of weights $\leq m-2n$. If $\mathcal X$ is smooth and proper over $\F_q$, then $H^m(\mathcal X_{\F},\Ql(n))$ is pure of weights $m-2n$.
\end{fact}
By (\ref{inclusionequation}) and Fact \ref{weights}, the cokernel of $\widetilde{\mathfrak i}^*:H^{2i}(\widetilde{\mathcal Z}^{cmp}_{\F},\Ql(i))\rightarrow H^{2i}(\mathcal Z_{\F},\Ql(i))$ is mixed of weights $\geq 1$, while $H^{2i}(\widetilde{\mathcal Z}^{cmp}_{\F},\Ql(i))$ is pure of weight $0$. Hence, the image of $\widetilde{\mathfrak i}^*:H^{2i}(\widetilde{\mathcal Z}^{cmp}_{\F},\Ql(i))\rightarrow H^{2i}(\mathcal Z_{\F},\Ql(i))$ consists exactly of the generalized eigenspace on which $\varphi$ acts with generalized eigenvalues of weight $0$. So $\mathfrak h^*(z')\in H^{2i}(\widetilde{\mathcal Z}_{\F},\Ql(i))^{\varphi}_{gen}$ is in the image of $\widetilde{\mathfrak i}^*:H^{2i}(\widetilde{\mathcal Z}^{cmp}_{\F},\Ql(i))\rightarrow H^{2i}(\mathcal Z_{\F},\Ql(i))$, hence by (\ref{usoprociclico}), $\mathfrak h^*(z')$ is the image of some $\widetilde z\in H^{2i}(\widetilde{\mathcal Z}^{cmp}_{\F},\Ql(i))^{\varphi}_{gen}$.
\subsubsection{Using the Tate conjecture}\label{usingtate}
Since $WS(\widetilde{\mathcal Z}^{cmp},i,\ell)$ holds by assumption, the injection
$$H^{2i}(\widetilde{\mathcal Z}^{cmp}_{\F},\Ql(i))^{\pi_1(\F_q)}\hookrightarrow H^{2i}(\widetilde{\mathcal Z}^{cmp}_{\F},\Ql(i))$$
has a $\pi_1(\F_q)$-equivariant splitting. So, since $\pi_1(\F_q)$ is procyclic generated by $\varphi$, one has
 $$H^{2i}(\widetilde{\mathcal Z}^{cmp}_{\F},\Ql(i))^{\varphi}_{gen}=H^{2i}(\widetilde{\mathcal Z}^{cmp}_{\F},\Ql(i))^{\varphi}=H^{2i}(\widetilde{\mathcal Z}^{cmp}_{\F},\Ql(i))^{\pi_1(\F_q)}.$$
Hence, by $T(\widetilde{\mathcal Z}^{cmp},i,\ell)$, there exists a $\widetilde w\in CH^i(\widetilde{\mathcal Z}_{\F}^{cmp})\otimes \Ql$ such that $c_{\widetilde{\mathcal Z}_{\F}^{cmp}}(\widetilde w)=\widetilde z$. We conclude the proof observing that, thanks to the commutative diagram at the end of \ref{sspread}, $h^*(z)$ is the image of $\widetilde i_{\eta}^*\widetilde{\mathfrak i}^*(\widetilde w)$ via $c_{\widetilde Z_{\overline k}}:CH^i(\widetilde Z_{\overline k})\otimes \Ql\rightarrow H^{2i}(\widetilde Z_{\overline k},\Ql(i)).$
This concludes the proof of Proposition \ref{propositionkey}.

\vskip 1\baselineskip
 \textit{emilianomaria.ambrosi@gmail.com}\\
Centre de Math\'ematiques Laurent Schwartz - Ecole Polytechnique,\\
91128 PALAISEAU, FRANCE.

\begin{thebibliography}{99}
\bibitem[Amb18]{mioneron} E. Ambrosi,
Specialization of N\'eron-Severi groups in positive characteristic, preprint available \href{http://emiliano.ambrosi.perso.math.cnrs.fr/}{here}.
\bibitem[And96a]{andk3} Y. Andr\'e, On the Shafarevich and Tate conjectures for hyperkähler varieties. Math.
Ann. 305 p. 205-248, 1996.
\bibitem[And96b]{AndreIHES} Y. Andr\'e, 
Pour une th\'eorie inconditionnelle des motifs, Publ. Math. IHES 83, p. 5-49, 1996.
\bibitem[And04]{andremotivi} Y. Andr\'e, Une introduction aux motifs (motifs purs, motifs mixtes, p\'eriodes), Panorama et synth`ese 17, S.M.F. 2004.
\bibitem[And06]{Andre} Y. Andr\'e, 
Deformation et specialisation de cycle motives, Journal of the Institute of Mathematics of Jussieu 4, p. 563-603, 2006.
\bibitem[Cad12]{MC} A. Cadoret, 
Motivated cycles under specialization, In "Groupes de Galois Géométriques et differentiels", P. Boalch and J.-M. Couveignes eds., S\'eminaires et Congrès, S.M.F., p. 25-55, 2012.
\bibitem[Cha13]{charlietatek3} F. Charles, The Tate conjecture for K3 surfaces over finite fields, Invent. Math. 194, p. 119-145, 2013.
\bibitem[dJ]{djtate} A.J. de Jong, Tate conjecture for divisors, unpublished note.
\bibitem[dJ96]{alteration} A.J. de Jong, Smoothness, semi-stability and alterations. Publ. Math. I.H.E.S. 83, p. 51–93, 1996.
\bibitem[Del68]{Spectraldeligne} Deligne, P. Theoreme de Lefschetz et criteres de degenerescence de suites spectrales. Publ. Math. IHES 35, p.107-126, 1968.
\bibitem[Del80]{Weil2} P. Deligne,  La Conjecture de Weil, II, Publ. Math. IHES 52, p. 137-252, 1980.
\bibitem[FW84]{faltingsabelian} G. Faltings, G. Wustholz (eds.), Rational Points, Aspects of Mathematics E6, Friedr. Vieweg, 1984.
\bibitem[Fu99]{leifusiccomeimmobile} L. Fu, On the semisimplicity of pure sheaves. Proc. Amer. Math. Soc. 127, p. 2529-2533, 1999.
\bibitem[KMP15]{peratatek32} W. Kim and K. Madapusi Pera, 2-adic integral canonical models and the Tate conjecture in characteristic 2, preprint, 2015.
\bibitem[LMS14]{finitenessk3} M. Lieblich, D. Maulik, A. Snowden, Finiteness of K3 surfaces and the Tate conjecture, Ann. Sci. Ec. Norm. Super. 2, p. 285-308, 2014.
\bibitem[Moo17]{Moonentate} B. Moonen, On the Tate and Mumford-Tate conjectures in codimension one for varieties with $h^{2,0}=1$, Duke Math J. 166, p. 739-799, 2017.
\bibitem[Moo18]{Moonenss} B.Monnen, A remark on the Tate conjecture, to appear in the J. of Algebraic Geometry.
\bibitem[Mor15]{Varmor}
M. Morrow,  A Variational Tate Conjecture in crystalline cohomology, to appear in Journal of the European Mathematical Society.
\bibitem[MP15]{peratatek3} K. Madapusi Pera, The Tate conjecture for K3 surfaces in odd characteristic, Invent. Math. 201, p. 625-668, 2015.
\bibitem[NO85]{Nygaardk3} N.O. Nygaard and A. Ogus, Tate’s conjecture for K3 surfaces of finite height, Ann. Math. 122, p. 461-507, 1985.
\bibitem[Ser89]{Serreweil}
J. P. Serre, Lectures on the Mordell-Weil theorem, Aspects of Mathematics, Friedr. Vieweg and Sohn, 1989.
\bibitem[SP]{stacks-project}
The Stacks Project Authors, Stacks Project, available \href{http://stacks.math.columbia.edu}{here}
\bibitem[Tan95]{tank3} S. G. Tankeev. Surfaces of type K3 over number fields and the Mumford–Tate conjecture. Izv. Ross. Akad. Nauk Ser. Mat. 59, 179-206, 1995.
\bibitem[Tat65]{tateconj} J. Tate, Algebraic cohomology classes and poles of zeta functions, Arithmetical Algebraic Geometry, New York: Harper and Row,  p. 93-110, 1965.
\bibitem[Tat66]{tateabelianfinite} J. Tate, Endomorphisms of abelian varieties over finite fields, Invent. Math. 2, p. 134-144, 1966
\bibitem[Tat94]{Tate} J. Tate, Conjectures on algebraic cycles in $\ell$-adic cohomology, Motives, Proc. Sympos. Pure Math., Part 1, Amer. Math. Soc., p. 71-83, 1994.
\bibitem[Zar75]{zarhinabelian1} Yu.G. Zarhin, Endomorphisms of abelian varieties over fields of finite characteristic, Izv. Akad. Nauk SSSR Ser. Mat. 39, p. 272-277, 1975.
\bibitem[Zar77]{zarhinabelian2} Yu. G. Zarhin, Endomorphisms of abelian varieties and points of finite order in characteristic p, Mat. Zametki 21, p. 737-744, 1977.
\end{thebibliography}
\end{document}